\documentclass[12pt]{article}

\usepackage{amssymb,amsmath}
\def\cc {{\mathfrak c}}
\def\fU {{\mathfrak U}} 
\def\TT {{\mathbb T}}
\def\gs {{\sigma}}

\input xy  \xyoption{all}

\newfont{\fn}{eufm10 scaled \magstep1}

\newfont{\bn}{msbm10 scaled \magstep1}

\newfont{\bs}{msbm7 scaled \magstep1}

\title{ Products of finite groups\\ and nonmeasurable subgroups}

\author{ F. Javier Trigos-Arrieta}

\begin{document}

\maketitle

\begin{abstract}
\noindent It is proven that if $G$ is a finite group, then $G^\omega$ has $2^\cc$ dense nonmeasurable 
subgroups. Also, other examples of compact groups
with dense nonmeasurable subgroups are presented.
\end{abstract}

\section{Introduction }

In \cite{SaStr1985}, the authors asked whether every infinite compact group has a (Haar) nonmeasurable (dense) subgroup. That every Abelian infinite compact group does is proven in \cite{HR} (16.13(d)). That every
non-metric compact group bigger than $\cc$ does follows from the fact that every such group has a proper pseudocompact subgroup \cite{Itz-Shakh-97}, which in turn is nonmeasurable \cite{Comfort84} (6.14). Thus, the problem remains open only for non-abelian metric and non-metric  groups of cardinality $\cc$. In this short note we prove the result in the abstract, and using \cite{CRT} (2.2) show that the unitary groups $\fU(n)$ do have too dense nonmeasurable subgroups. 

\section{Unitary groups }

The result \cite{CRT} (2.2) states that if $K$ and $M$ are compact groups and $\varphi:K \to M$ is a continuous homomorphism onto, then the preimage of any (dense)  nonmeasurable subgroup of $M$ is a (dense)  nonmeasurable subgroup of $K$. Since the torus $\TT$ has plenty of (dense)  nonmeasurable subgroups, and the determinant is a continuous homomorphism from any unitary group $\fU(n)$ \cite{HR} (2.7(b)) onto $\TT$, it follows that the unitary groups do have dense nonmeasurable subgroups. 

\newpage

\section{Countable products of finite groups }

\def\sU {{\mathcal U}}
\def\sI {{\mathcal I}}
Let $\sU$ be a free ultrafilter. Consider $\sI:=2^\omega \setminus \sU$. The collection $\sI$ will be called an {\em ideal}. The following are properties dual of those for an ultrafilter:
\begin{enumerate}
\item $A \subset \omega \implies \omega\setminus A \in \sI,$ or $A \in \sI$,
\item $A \in \sI \implies \omega\setminus A \not\in \sI,$
\item$A \in \sI, C \subseteq A \implies C \in \sI,$ and
\item $A, B \in \sI \implies A \cup B \in \sI.$
\end{enumerate}

For each $n \in \omega$, let $G_n$ be a non-trivial finite group, with identity $e_n$. Consider $G:=\times_{n < \omega} G_n$. If $x=(x_n) \in G$, denote by $\gs(x):=\{n < \omega:x_n\neq e_n\}$. If $A \subseteq \omega$, let $G_A:=\{x\in G:\gs(x) \subseteq A\}$. 
Finally, denote by $G_\sI:=\cup_{A\in \sI} G_A$. 
Clearly, $G_\sI$ is a subgroup of $G$, and because $\sU$ is a free ultrafilter, $G_\sI$ is dense in $G$.

\noindent {\bf Question (3.1)} {\em Is $G_\sI$ a measurable subgroup of $G$?}
\def\st {such that}
\def\ie {i. e.}
\def\nhd {neighbourhood}
\def\pf {{\em \noindent Proof:}} 
\def\qed {{\blacksquare}}
\def\epf {~\hfill $\qed$} 

We can answer this question, negatively, if all $G_n$ are equal, say to $\Gamma$. Denote by $e$ the identity of $\Gamma$. First of all, we will prove that, in this case, $G/G_\sI\simeq\Gamma$. Let $x \in G$. For each $a\in \Gamma\setminus \{e\}$, denote by $\gs(x,a)$ those $n\in\gs(x)$ \st \ $x_n=a$. Notice therefore that $\gs(x)$ is the disjoint union of the $\gs(x,a)$ as $a$ runs through every non-identity element in $\Gamma$.

If $x\not \in G_\sI$, then $\gs(x) \not\in \sI$. We claim that there is a unique $a\in\Gamma\setminus \{e\}$ with $\gs(x,a)\not\in \sI$. For, if for each $a\in \Gamma\setminus \{e\}$, we had that $\gs(x,a)\in \sI$, then we would have $\gs(x)\in \sI$, 
a contradiction. Thus there is $a_0\in \Gamma\setminus \{e\}$ with $\gs(x,a_0)\not\in \sI$. Hence $\omega \setminus \gs(x,a_0)\in \sI$, and since $a\in \Gamma\setminus \{e,a_0\} \implies \gs(x,a) \subseteq \omega \setminus \gs(x,a_0)$, the properties for ideals show that 
$\cup_{a\in \Gamma\setminus \{e,a_0\}}\gs(x,a) \in \sI$. Now, define $y=(y_n)$ by
\[y_n:=\left\{ \begin{array}{lr}
 a_0^{-1}x_n,\: \mbox{if}\: n\in \cup_{a\in \Gamma\setminus 
			\{e,a_0\}}\gs(x,a), \\
e,\: \mbox{if}\: n \in \gs(x,a_0), \\
 a_0^{-1},\: \mbox{otherwise.} 
\end{array} \right. \] 

Because, $\gs(y)=\omega \setminus \gs(x,a_0) \in \sI$, it follows that $y \in G_\sI$. Set $\overline{a_0}=(t_n)$ by $t_n:=a_0$ for all $n<\omega$, i.e., it's the constant sequence $a_0$. We now show that 
\[x=\overline{a_0}\cdot y.\]
For, if $n\in \cup_{a\in \Gamma\setminus 
			\{e,a_0\}}\gs(x,a)$, then $t_n a_0^{-1}x_n=a_0a_0^{-1}x_n=x_n$. If $n \in \gs(x,a_0)$, then $t_n e=a_0e=a_0$. And if $n \not\in \cup_{a\in \Gamma\setminus 
			\{e\}}\gs(x,a)$, then $t_n a_0^{-1}=a_0a_0^{-1}=e=x_n$, as required.

This shows the following:

\noindent {\bf Theorem (3.2)} {\em If $\Gamma$ is a finite group, and $G:=\Gamma^\omega$, then $G/G_\sI\simeq\Gamma$.}\epf

Thus $G_\sI$ has finite index and therefore cannot have zero measure. 

\noindent {\bf Theorem (3.3)} {\em (Steinhaus-Weil Theorem) If $F$ is a measurable subset of a (locally) compact group $G$ with strictly positive (left Haar) measure, then $F \cdot F^{-1}:=\{x y^{-1}: x, y \in F\}$ contains a \nhd \ of the identity of $G$. Thus, if $F$ is in addition a dense subgroup of $G$, then $F=G$. \epf}

This is proven in \cite{weil40}. See also \cite{steinhaus20} and \cite{strom72}.

\noindent {\bf Corollary (3.4).} {\em $G_\sI$ is not measurable.}

\pf \ If $G_\sI$ were measurable, then it would have strictly positive measure. By the above theorem,  it would have to be equal to the whole $G$, clearly a contradiction. \epf
\def\AA {{\mathbb A}}

Now, assume that $\Gamma$ is a simple (finite) non-Abelian group (for example, the alternating subgroup $\AA_m$  on $m$ elements, with $m\geq 5$). Robert Bassett and the author have proved that the only normal subgroups of $G$ are of the form $G_\sI$ for some ideal $\sI$. If we continue assuming that $\sI$ is the complement in $2^\omega$ of a free ultrafilter, then it follows that $G_\sI$ is a {\em maximal normal} subgroup. Let $\varphi: G \to G/G_\sI$ be the natural map. Identify, by Theorem 1, $G/G_\sI$ with $\Gamma$ and $G_\sI$ with $e$. Choose $g\in \Gamma, g \neq e$ and denote by $B$ the subgroup of $\Gamma$ generated by $g$. Because $\Gamma$ is simple and non-Abelian, $\{e\} \subset B \subset \Gamma$ and these contentions are proper. Set $H:=\varphi^{\leftarrow}[B] $. Thus $G_\sI \subset H \subset G$, with the above contentions proper. Since $G_\sI$ is a maximal {\em normal} subgroup properly contained in $H$, it follows that $H$ is a non-normal subgroup of $G$. And since the contention $H \subset G$ is proper, it follows that $H$ is a non-normal {\em proper} subgroup of $G$. Hence, another application of Steinhaus'-Weil Theorem (3.3) implies the following:

\noindent {\bf Corollary (3.5)} {\em $H$ is a non-normal not measurable subgroup of $G$.} \epf

\def\< {{\langle}}
\def\> {{\rangle}}
\noindent {\bf Example (3.6)} The condition that all $G_n$ are equal in Corollary (3.4) is necessary as this example shows. Let $ \< t_n \> _{n< \omega}$ be a an increasing sequence of non-zero numbers converging to 1, such that $g_m:=t_0\cdot t_1\cdot t_2 \cdots t_{m-1}$ converges to say $t \in (0,1)$ (for example, if $\Sigma_{n=0}^\infty a_n$ converges with $1> a_n \downarrow 0$, then $t_n:=1-a_n$ satisfies the condition, see Stromberg's book \cite{strom81}). Now, pick a strictly increasing sequence of integers $ \< k_n \> _{n< \omega} $ \st \ $t_n \leq \frac{k_n-1}{k_n} < 1$. If $\tau_n:=\frac{k_n-1}{k_n}$, then $\gamma_m:=\tau_0\cdot \tau_1 \cdots \tau_{m-1}$ converges to say $\gamma \in [t,1)$. Set $G_n:=\AA_{k_n}$, and of course $G:=\times_{n<\omega} G_n$. Denote by $m$ the (Haar) measure on $G$. We claim that $m(G_\sI )=0$.  To see this, denote by $1_n$ the identity of $G_n$. Set $\omega(n):=\omega \setminus n=\{n, n+1,...\}$, and $B_n:=\{x\in G:\omega(n)\subseteq \gs(x)\}$. Basically, $B_n$ consists of those $x$ whose first $n$ coordinates can be anything, but everything after must be different than $1_n$. Notice then that $B_n=G_0 \times G_1 \times \cdots \times G_{n-1} \times (\times_{k\geq n} (G_k \setminus \{1_k\}))$, hence $B_0 \subseteq B_1 \subseteq \cdots $, and therefore, $G\setminus B_0 \supseteq G\setminus B_1 \supseteq \cdots $. Since the measure of $G_n \setminus \{1_n\}$, in $G_n$, is $\frac{k_n-1}{k_n}$, it follows that $m(B_n)=\lim_{m\to \infty} \Pi_{j=n}^{m-1}\frac{k_j-1}{k_j}=(\frac{\tau_0\cdot \tau_1 \cdots \tau_{n-1}}{\tau_0\cdot \tau_1 \cdots \tau_{n-1}})(\lim_{m\to \infty}(\tau_n\cdots \tau_{m-1}))=(\frac{1}{\gamma_n})(\lim_{m\to \infty}(\tau_0\cdot \tau_1 \cdots \tau_{n-1}\tau_n\cdots \tau_{m-1}))=\frac{\gamma}{\gamma_n}$. Thus $m(G\setminus B_n)=1-\frac{\gamma}{\gamma_n}$. Since $\omega(n) \not \in \sI$, for all $n< \omega$, we have that $G_\sI \subseteq \bigcap_{n<\omega} (G\setminus B_n)$, which, by Proposition 2, Chapter 11 in \cite{Royden1968}, has measure $\lim_{n<\omega} m(G\setminus B_n)=\lim_{n<\omega} (1-\frac{\gamma}{\gamma_n})=1-\frac{\gamma}{\lim_{n<\omega}\gamma_n}=1-\frac{\gamma}{\gamma}=0$. Therefore $G_\sI$, in this case, has measure 0, as required. 

Nevertheless, Corollary (3.4) can be improved as follows, by using \cite{CRT} (2.2):

\noindent {\bf Corollary (3.7)} {\em For each $n \in \omega$, let $G_n$ be a non-trivial finite group \st \ $G_n=\Gamma$, some fixed group $\Gamma$, for infinitely many $n \in \omega$. Then $G:=\times_{n < \omega} G_n$ has nonmeasurable subgroups.}

\pf \ Let $\omega_\Gamma:=\{n \in \omega: G_n =\Gamma\}$. By Corollary 1, $\Gamma^\omega$ has nonmeasurable subgroups, and since $G_\Gamma:=\times_{n \in \omega_\Gamma} G_n$ is topologically isomorphic to $\Gamma^\omega$, it does 
too have  nonmeasurable subgroups. Since $G=\times_{n < \omega} G_n = G_\Gamma \times (\times_{n < \omega\setminus \omega_\Gamma} G_n)$ , the
projection of $G$ onto the first factor, yields the result. \epf

\section{Final Remarks }
\begin{enumerate}
\item That unitary groups have nonmeasurable subgroups was obtained during a wonderful dinner in Middletown back in 2002, when the author met with his teachers and friends, Wis Comfort, Tony Hager and Lew Robertson.

\item Faculty in the Department of Mathematics at CSUB made the author aware of a mistake in an older version of Example 1. 

\item S. Hern\' andez has communicated to the author that he, K. Hofmann and S. Morris have independently generalized most of the results in this article,  with quite different techniques.
\end{enumerate}

\noindent{\small\em Department of Mathematics}\\
{\small\em California State University, Bakersfield}\\
{\small\em Bakersfield, California, USA}\\
{\small\em e-mail:  jtrigos@csub.edu}
\end{document}